\documentclass[11pt]{amsart}

\usepackage{amsfonts, amsmath, amssymb}

\newcommand{\E}{\mathsf{E}}
\newcommand{\Pp}{\mathsf{P}}
\newcommand{\R}{\mathbb{R}}
\newcommand{\N}{\mathbb{N}}
\newcommand{\Z}{\mathbb{Z}}
\newcommand{\cov}{\mathop{\mathsf{cov}}}
\newcommand{\supp}{\mathop{\mathrm{supp}}}
\newcommand{\ONE}{{\bf 1}}

\newtheorem{theorem}{Theorem}
\newtheorem{remark}{Remark}

\title{Poisson Limit for Associated Random Fields}
\author{Yuri Bakhtin}

\begin{document}
\maketitle
\begin{abstract}We prove that under an easily verifiable set of conditions a sequence of associated random fields converges
under rescaling to the Poisson Point Process and give a couple of examples.
\end{abstract}

\section{Introduction And Main result}

In this note we prove a Poisson scaling limit for a sequence of associated random fields.
Let us recall that a finite family (vector) $(X_1,\ldots,X_m)$ of random variables (r.v.'s) is called
associated if for every pair of bounded and coordinatewise nondecreasing
functions $f,g:\R^m\to\R$,
\begin{equation}
 \cov(f(X_1,\ldots,X_m),g(X_1,\ldots,X_m))\ge 0.
\label{eq:fkg}
\end{equation}
An infinite family of r.v.'s is called associated if its every finite subfamily is associated.

The notion of association was introduced and studied in~\cite{EPW-MR0217826}.  Inequalities~\eqref{eq:fkg} with their equivalents
have been often referred to as FKG inequalities by the initials of authors of~\cite{FKG-MR0309498} who studied this type of positive correlation independently.

Associated r.v.'s arise
frequently in various problems of statistical mechanics and many other areas, see numerous examples,
a historic overview, theory and applications in a recent monograph~\cite{Bu-Sh-MR2375106}.

Basic properties of associated random vectors: jointly independent r.v's form an associated family; monotone transformations of
associated random vectors are associated, too.

A number of limit theorems for sums of associated r.v.'s have been proved, see~\cite{Bu-Sh-MR2375106} and references therein. To the best of our knowledge, no theorem
on convergence to a Poisson Point Process has appeared in the literature, although some results on Poisson approximations
for systems satisfying FKG inequalities can be found in \cite{Ganesh-et-al:MR1762669} and references therein.

We proceed to describe the setting. We fix a dimension $d\in\N$, and for each $n\in\N$, let $(X^{(n)}_j)_{j\in \Z^d}$ be a weakly stationary (i.e. in the sense of first
moment and covariance)
associated random field. We assume that
for all $n,j$, r.v.~$X^{(n)}_j$ takes two values, 0 and 1,
 and there is a number $\lambda>0$
such that
\[
 p_n=\frac{\lambda+o(1)}{n^d},
\]
where $p_n=\Pp\{X^{(n)}_0=1\}$.

We also assume that
\begin{equation}
\label{eq:decay_of_sigma}
 \lim_{n\to\infty} n^d \sigma(n)=0,
\end{equation}
where
\[
\sigma(n)=\sum_{j\ne 0}\cov(X_0^{(n)},X_j^{(n)}).
\]

For any $n$ we define a random measure $\mu_n$ on $\R^d$ via
\[
 \mu_n(A)=\sum_{j\in \Z^d \cap nA} X^{(n)}_j,
\]
where $nA=\{nx:\ x\in A\}$.

The vague topology on locally bounded Borel measures is defined by its base, the class of finite intersections of sets of
the form $\{\nu: s<\int_{\R^d}fd\nu<t \}$ with arbitrary nonnegative continuous function $f$ with bounded support and $s,t\in\R$,
see~\cite[Appendix 7]{Kallenberg-MR0431373}.

\begin{theorem}\label{th:conv_to_poisson}
 Under the conditions stated above, the sequence of measures~$\mu_n$ converges in distribution 
in the vague topology to the Poisson measure $\mu$ with parameter $\lambda$.
\end{theorem}

\begin{proof}
By \cite[Theorem 4.2]{Kallenberg-MR0431373}, it is sufficient to check that
for every continuous
nonnegative function $f$ with compact support,
\[
\int_{\R^d} fd\mu_n\stackrel{Law}{\to} \int_{\R^d} fd\mu,\quad\mbox{as}\ n\to\infty.
\]

Take a continuous function $f$ with compact support and a number $t\in\R$, and find
\begin{align*}
 \E e^{it\int fd\mu_n}&=\E e^{it\sum_{j\in\Z^d} f(\frac{j}{n})X^{(n)}_j}\\
&=\prod_{j\in\Z^d} \E e^{it f(\frac{j}{n})X^{(n)}_j}+
\left|\E e^{it\sum_{j\in\Z^d} f(\frac{j}{n})X^{(n)}_j}-\prod_{j\in\Z^d} \E e^{it f(\frac{j}{n})X^{(n)}_j} \right|\\
&=I_1(n)+I_2(n).
\end{align*}

Notice that, in fact, the product in $I_1(n)$ involves finitely many factors, and
\[
 I_1(n)=\prod_{j\in\Z^d}\left(1+p_n(e^{itf(\frac{j}{n})}-1)\right).
\]
Choosing the main branch of the natural logarithm $\ln$, we can write
\[
 I_1(n)=\exp\left\{\sum_{j\in\Z^d}\ln(1+p_n(e^{itf(\frac{j}{n})}-1))\right\}.
\]
Using the boundedness of $f$ and the Taylor expansion for the logarithm we derive that
\[
 I_1(n)=\exp\left\{\frac{\lambda+o(1)}{n^d}\sum_{j\in\Z^d}(e^{itf(\frac{j}{n})}-1)\right\}(1+o(1)).
\]
Obviously, the r.h.s converges to 
\[
 \phi(t)=\exp\left\{\lambda\int_{\R^d}(e^{itf(x)}-1)dx\right\},
\]
the characteristic function of $\int_{\R^d} fd\mu$, and the proof will be finished as soon as we show that 
\begin{equation}
\lim_{n\to\infty} I_2(n)=0.
\label{eq:lim_I_2}
\end{equation}

To estimate $I_2(n)$ we need Newman's inequality:
\begin{theorem}[\cite{Newman-MR576267}]\label{th:newman_inequality}
If $(Y_1,\ldots,Y_m)$ is a family of associated r.v.'s with finite second moment then
\[
\left|\E e^{i\sum_{j=1}^m r_j Y_j} - \prod_{j=1}^m \E e^{ir_j Y_j}\right|\le\frac{1}{2}\sum_{j_1\ne j_2}|r_{j_1}r_{j_2}|\cov(Y_{j_1},Y_{j_2}),
\]
for any real numbers $r_1,\ldots, r_m$.
\end{theorem}
Applying this inequality to $I_2(n)$ we see that
\begin{align*}
 I_2(n)&\le \frac{t^2 \|f\|^2_{L^\infty}}{2}\sum_{\substack{ j_1,j_2\in\Z^d\cap n\supp(f)\\j_1\ne j_2} }\cov(X^{(n)}_{j_1},X^{(n)}_{j_2}),\\
&\le \frac{t^2 \|f\|^2_{L^\infty}}{2}\, |\Z^d\cap n\supp(f)|\, \sigma(n),
\end{align*}
where $\supp(f)$ denotes the support of $f$, and $|\cdot|$ denotes the number of elements. Since
$|\Z^d\cap n\supp(f)|\le K n^d$ for some constant $K>0$ and all $n>0$, \eqref{eq:lim_I_2} follows
from \eqref{eq:decay_of_sigma}.
\end{proof}

\begin{remark} The crucial step in the proof above is the application of Newman's inequality for associated
random variables. Covariance inequalities of this type can be obtained for a wide class of dependent r.v.'s.
In particular the theorem is also applicable if one replaces association by quasi-association, 
see~\cite{Bu-Su-MR1858636} and proof of Theorem~\ref{th:newman_inequality} in \cite{Newman-MR576267}.
\end{remark}

\section{Examples}
Let $G$ be a finite subset of $\Z^d$ for some $d\in\N$. Denote $m=|G|$ and for each $n$ consider an i.i.d. family
$(Y^{(n)}_k)_{k\in\Z^d}$ of Bernoulli random variables with 
\[
\Pp\{Y^{(n)}_0=x\}=\begin{cases}
                    q_n,& x=1,\\
                    1-q_n,& x=0,
                   \end{cases}
\]
where 
\[
q_n=\frac{1}{n^{d/m}}.
\]

For any finite subset $H$ of $\Z^d$ and every $n$, we denote 
\[
\chi^{(n)}_H=\prod_{j\in H}Y^{(n)}_j=\ONE_{\{Y_j^{(n)}=1,\ j\in H\}},
\]
and define a random field $(X^{(n)}_k)_{k\in\Z^d}$ via
\[
 X^{(n)}_k= \chi^{(n)}_{k+G},
\]
where $k+G=\{k+j:\ j\in G\}$. Poisson approximations for a similar model with rectangular $G$ was considered in \cite{Fu-Koutras:MR1272757}.

Let us verify that $X^{(n)}$ satisfies the conditions of Theorem~\ref{th:conv_to_poisson}.
Random field $Y^{(n)}$ is associated since it is composed of independent components. 
Therefore, $X^{(n)}$ is associated being a monotone transform of the
associated field $Y^{(n)}$. It is also stationary due to stationarity of $Y^{(n)}$.

For each $n$, $X^{(n)}_0$ is a Bernoulli r.v. with
\[
\Pp\{X^{(n)}_0=1\}= \Pp\{Y_j^{(n)}=1,\ j\in G\}=\left(\frac{1}{n^{d/m}}\right)^m=\frac{1}{n^d}.
\]
Let us now estimate $\sigma(n)$. Notice that $\cov(X^{(n)}_0,X^{(n)}_j)=0$ for sufficiently large values of $|j|$,
so that there is a number $M$ such that for all $n$,
\begin{equation}
 \sigma(n)\le M \max_{j\ne 0}\cov(X^{(n)}_0,X^{(n)}_j).
\label{eq:bound_on_sigma}
\end{equation}

Notice that
\[
 \cov(X^{(n)}_0,X^{(n)}_j)= \E \chi^{(n)}_{G\cup (j+G)}-\E \chi^{(n)}_{G}\E \chi^{(n)}_{j+G}.
\]
Since a finite set cannot be invariant under a translation, $|G\cup (j+G)|\ge m+1$ for any $j$. Therefore,
\[
 \cov(X^{(n)}_0,X^{(n)}_j)\le \frac{1}{n^{d(m+1)/m}}=o(1/n^d),
\]
which, together with \eqref{eq:bound_on_sigma}, implies \eqref{eq:decay_of_sigma}, so that all the
conditions of Theorem~\ref{th:conv_to_poisson} are satisfied.

For an associated random field $X^{(n)}$, condition~\eqref{eq:decay_of_sigma} means that $X^{(n)}_0$ is asymptotically independent of the rest of the random field. There is a variety of situations that can
happen if this condition is replaced with weaker restrictions on dependence. The next example illustrates
the convergence to a compound  Poisson point process (with nonrandom mass 2 assigned to each atom), see~\cite{DV-MR950166} for the definition and properties of compound Poisson point processes.

Consider $d=1$, and for every $n$ and all $k\in\Z$,
\[
X^{(n)}_{k}=Y^{(n)}_{k}\vee Y^{(n)}_{k+1}=Y^{(n)}_{k}+Y^{(n)}_{k+1}-Y^{(n)}_{k}Y^{(n)}_{k+1},
\] 
where $Y^{(n)}$ is a sequence of i.i.d. Bernoulli r.v.'s with
$\Pp\{Y^{(n)}_{0}=1\}=1/n$. Then, as an easy computation shows, $\sigma(n)\sim 1/n$ so that \eqref{eq:decay_of_sigma} is violated. One can also show that the sequence of random measures
$\mu_n$ converges in distribution to $2\mu$, where $\mu$ is the Poisson process with unit intensity, so that
the conclusion of Theorem 1 is violated as well. Indeed, take a  continuous function $f$ with compact support, and write
\[
 \E e^{it\int_{\R}fd\mu_n}=\E e^{it\sum_{j\in\Z} f(\frac{j}{n})(Y^{(n)}_{j}+Y^{(n)}_{j+1})-it\sum_{j\in\Z} f(\frac{j}{n})Y^{(n)}_{j}Y^{(n)}_{j+1}}.
\]
Notice that
\[
\sum_{j\in\Z} f\left(\frac{j}{n}\right)Y^{(n)}_{j}Y^{(n)}_{j+1}\stackrel{\Pp}{\to}0,\quad n\to\infty,
\]
due to the Markov inequality, since the expectation of l.h.s. is $O(1/n)$. Therefore, we see that
 \begin{align*}
 \lim_{n\to\infty}\E e^{it\int_{\R}fd\mu_n}&=\lim_{n\to\infty}\E e^{it\sum_{j\in\Z} (f(\frac{j}{n})+f(\frac{j-1}{n}))Y^{(n)}_{j}}\\
&=\exp\left\{\int_{\R}(e^{it2f(x)}-1)dx\right\},
\end{align*}
by the same argument we used to analyze $I_1(n)$.
\bibliographystyle{plain}
\bibliography{fkg-poisson}

\begin{thebibliography}{1}

\bibitem{Bu-Sh-MR2375106}
Alexander Bulinski and Alexey Shashkin.
\newblock {\em Limit theorems for associated random fields and related
  systems}.
\newblock Advanced Series on Statistical Science \& Applied Probability, 10.
  World Scientific Publishing Co. Pte. Ltd., Hackensack, NJ, 2007.

\bibitem{Bu-Su-MR1858636}
Alexander Bulinski and Charles Suquet.
\newblock Normal approximation for quasi-associated random fields.
\newblock {\em Statist. Probab. Lett.}, 54(2):215--226, 2001.

\bibitem{DV-MR950166}
D.~J. Daley and D.~Vere-Jones.
\newblock {\em An introduction to the theory of point processes}.
\newblock Springer Series in Statistics. Springer-Verlag, New York, 1988.

\bibitem{EPW-MR0217826}
J.~D. Esary, F.~Proschan, and D.~W. Walkup.
\newblock Association of random variables, with applications.
\newblock {\em Ann. Math. Statist.}, 38:1466--1474, 1967.

\bibitem{FKG-MR0309498}
C.~M. Fortuin, P.~W. Kasteleyn, and J.~Ginibre.
\newblock Correlation inequalities on some partially ordered sets.
\newblock {\em Comm. Math. Phys.}, 22:89--103, 1971.

\bibitem{Fu-Koutras:MR1272757}
James~C. Fu and Markos~V. Koutras.
\newblock Poisson approximations for {$2$}-dimensional patterns.
\newblock {\em Ann. Inst. Statist. Math.}, 46(1):179--192, 1994.

\bibitem{Ganesh-et-al:MR1762669}
A.~Ganesh, B.~M. Hambly, Neil O'Connell, Dudley Stark, and P.~J. Upton.
\newblock Poissonian behavior of {I}sing spin systems in an external field.
\newblock {\em J. Statist. Phys.}, 99(1-2):613--626, 2000.

\bibitem{Kallenberg-MR0431373}
Olav Kallenberg.
\newblock {\em Random measures}.
\newblock Akademie-Verlag, Berlin, 1976.

\bibitem{Newman-MR576267}
C.~M. Newman.
\newblock Normal fluctuations and the {FKG} inequalities.
\newblock {\em Comm. Math. Phys.}, 74(2):119--128, 1980.

\end{thebibliography}

\end{document}